\documentclass[11pt]{article}
\usepackage{a4wide}
\usepackage[margin=1in]{geometry}
\usepackage[utf8]{inputenc}
\usepackage[english]{babel}
\usepackage{comment}
\usepackage{amssymb, amsmath}
\usepackage[colorlinks=true]{hyperref}
\usepackage{bbm}
\usepackage{mathtools}
\usepackage{latexsym}
\usepackage{amsfonts}
\usepackage{epsf}
\usepackage{psfrag}
\usepackage{epsfig}
\usepackage{textcomp}
\usepackage{color} 
\usepackage[normalem]{ulem}
\usepackage{enumerate}
\usepackage{tkz-euclide}

\newcommand{\R}{{\mathbb{R}}}   
\newcommand{\N}{{\mathbb{N}}}   
\newcommand{\E}{{\mathbb{E}}}   
\newcommand{\Prob}{{\mathbb{P}}}   
\renewcommand{\geq}{\geqslant}
\renewcommand{\leq}{\leqslant}

\usepackage{pgf}
\usetikzlibrary{arrows,automata}
\usepackage{subfig}
\usepackage{centernot}
\usepackage{mathrsfs,amsthm,amsbsy}

\usepackage[T1]{fontenc}
\usepackage{lmodern}
\theoremstyle{plain} 
\newtheorem{thm}{Theorem}[section] 

\theoremstyle{definition}

\newtheorem{lem}{Lemma}[section]

\theoremstyle{remark} 
\newtheorem{obs}{Remark}

\begin{document}

\title{Clustering in a preferential attachment network with triangles}
       
 \author{Angelica Pachon\thanks{Department of Mathematics, The Computational Foundry, Swansea University Bay Campus, Fabian Way, Swansea, SA1 8EN, United Kingdom, \url{a.y.pachon@swansea.ac.uk}}, \hspace{0.2cm} Robin Stephenson\thanks{School of Mathematical and Physical Sciences, University of Sheffield, Hicks Building, Hounsfield Road, Sheffield S3 7RH,
		United Kingdom, \url{robin.stephenson@normalesup.org}}}


\maketitle

  \begin{abstract}
 We study a generalization of the affine preferential attachment model where triangles are randomly added to the graph. We show that the model exhibits an asymptotically power-law degree distribution with adjustable parameter $\gamma\in (1,\infty)$, and positive clustering. However, the clustering behaviour depends on how it is measured. With high probability, the average local clustering coefficient remains positive, independently of $\gamma$, whereas the expectation of the global clustering coefficient does not vanish only when $\gamma>3$.  
 \end{abstract}

\textbf{keywords:} Random graph; Preferential attachment; Clustering coefficient

\textbf{MSC2020 subject classifications:}Primary: 05C80; 05C82. Secondary: 68R10; 91D30.

\section{Introduction}

The study of complex networks from a mathematical point of view uses random graphs and focuses on the most common structural properties, namely, the small-world nature and degree-distribution. The former is characterized by a high clustering and low distances, while the latter is usually a heavy-tailed distribution such as a power law, i.e. the fraction of vertices in the network with degree $k$, is close to $k^{-\gamma}$ for large values of $k$, with typically $2<\gamma<3$. 

In the area of network modelling using random graphs, a simple mechanism known as preferential attachment, was the first one to be used to model a network with an asymptotically power-law degree distribution with parameter $\gamma=3$. The model was proposed by Barab\'asi and Albert in \cite{BarabasiAlbert99} and later rigorously studied  in \cite{Bollobas2001,Cooper2003}.   This model is related to other previous models that appeared in the mathematical biology literature  \cite{Price,SIMON1955,Yule1925}. The preferential attachment model in \cite{BarabasiAlbert99} is a random graph process where at each time step, a new vertex is added and  connected to a given number of existing vertices with a probability proportional to their degrees. In \cite{Athreya} is  studied a generalization of this model by considering a probability proportional to a weight function of their degrees, $w(\cdot)$, so the model in  \cite{BarabasiAlbert99} corresponds to $w(x)=x$. The results in \cite{Athreya} indicate a different behaviour depending
on whether $w(\cdot)$ is  super-linear, sub-linear, or linear, and that the degree distribution is asymptotically a power law essentially only when the weight function is
linear: $w(x)=cx+\delta,$ $c>0$ and $c>-\delta$, where $\delta$ is a parameter for the initial \textit{attractiveness}. When the weight function is super-linear,  for example  $w(x) = cx^2$, the limit of the degree of any vertex as time goes to infinity exists and
is finite with probability 1, moreover, the empirical degree distribution  converges to the Dirac measure at 1.  When the weight function is sub-linear with $w(x)\sim cx^p$, $1/2<p<1$, then the degree of any vertex grows like $(\ln n)^q$, where $q=(1-p)^{-1}$ and hence there is no power-law growth.
This model has also been studied  in \cite{KaprivskyRedner01}. For  the super-linear and the sub-linear cases see also \cite{OliveiraSpencer05} and \cite{Rudas2004}, respectively, while for the linear case with $c=1$, see  \cite{BUCKLEY200453,Dorogovtsev2000,Mori03}, where it is showed that the asymptotic degree distribution follows a power law with parameter $\gamma=3+\delta$.
  
As mentioned before, another important property of real-world networks is \textit{clustering}, the tendency of a network to have groups of  densely interconnected vertices. There are several definitions for a clustering coefficient, the two most studied are the 'global clustering coefficient' and the 'average local clustering coefficient', which both measure the proportion of triangles present in a graph.  Creating random graph models that not only accurately capture a power law degree distribution, but also have a positive clustering coefficient (of either kind), has been a challenge. Notably, for both the model in \cite{Bollobas2001} and the model in \cite{Mori03} with initial positive attractiveness, the expected value of the global clustering coefficient is, up to a logarithmic correction, asymptotically equivalent to $\frac{1}{n}$ and therefore tends to $0$ (see respectively \cite{BollobasRiordan02} and  \cite{EGGEMANN2011953}).

In~\cite{HolmeKim02}, a variant of the model in \cite{Bollobas2001} was designed to exhibit clustering. The model essentially adds triangles to the graph. A rigorous analysis of this variant was recently carried out in \cite{Oliveira}, and it was shown that the average local clustering coefficient typically remains positive whereas the global clustering coefficient tends to zero at a slow rate. The authors also showed that this model satisfies an asymptotically power-law degree distribution with parameter $\gamma=3$. 

In the computer science literature, a common framework for the analysis of a wide class
of preferential attachment models was proposed in \cite{Ostroumova}. This framework includes the models in \cite{Bollobas2001}, \cite{HolmeKim02} and other generalizations.  It is constructed on the base of two properties. The first one is a restriction on the distribution of the increase of the degree of vertices at each time step, which is needed to get an asymptotically power-law degree distribution. The second one states that the graph should satisfies that the number of triangles should be proportional to the number of vertices. The results in \cite{Ostroumova} are non-rigorous but show that models which fit in this framework satisfy an asymptotically power-law degree distribution, and also give us the existence of three different regimes for the global clustering coefficient: $\gamma<3$, $\gamma=3$ and $\gamma>3$, where $\gamma$ is the parameter of the power-law decay.
 In the first two regimes the global clustering coefficient tends to zero, while when $\gamma>3$ it does not vanish as the graph grows. The authors also describe a specific model from this framework and use simulations to study its average local clustering coefficient. They show that it does not vanish in any of these three regimes. 

Taking into account all the previous results and ideas, in this paper we propose a simple modification to the preferential attachment model in \cite{Bollobas2001}, and rigorously study its degree distribution and clustering. Our model has two parameters, an initial attractiveness, $\delta>-1$, and a fixed probability to introduce triangles at each time step, $\alpha$.  We show that our the degree distribution of this model is asymptotically a power law  with parameter $\gamma=1+1/A$, where $A$ is a positive function which depends on $\delta$ and $\alpha$, and determines the three regimes, $\gamma<3$, $\gamma=3$ and $\gamma>3$, which correspond to $A>1/2$, $A=1/2$ and $A<1/2$, respectively. We prove that the expected value of the global clustering coefficient is of order $\Theta(1)$ if $A<1/2$, $\Theta(\ln(t)^{-1})$ if $A=1/2$ and $\Theta(t^{1-2A})$ if $A>1/2$, so in the limit it is only positive in the regime $\gamma>3,$ and is zero in the other two regimes. We also study the average local clustering coefficient and prove that, with high probability it does not vanish in any of the three regimes.  
 
The remainder of the paper is organized as follows. In Section 2, we describe the model. In Section 3, we present our main results and in Section 4 we prove these results. Section 5 concludes the paper.

\section{The Model}\label{TheModel}

The model that we will consider throughout this paper is a process of undirected simple graphs that we will call $(G(t),t\in\N_0).$ For $t\in\N_0$, $V_t$ and $E_t$ will respectively denote the vertex set and edge set of $G(t),$ with respective cardinalities $|V_t|$ and $|E_t|$. For $u$ and $v$ in $V_t$ for any $t$, we write $u\sim v$ if they are connected in $G(t),$ i.e. if $\{u,v\}\in E_t$ (as we will see, edges do not form between pre-existing vertices so the lack of mention of $t$ in the $\sim$ symbol is not an issue).


Our process has two parameters: a number $\alpha\in[0,1]$ which is the probability that we form a triangle at any step, and a number $\delta\in (-1,+\infty)$ which is the attractiveness of the vertices when not adding a triangle. We start at time $t=0$ with two vertices connected by an edge, that is $V_0:=\{v_1,v_2\}$ and $E_0:=\{(v_1,v_2)\}$. For $t\geq 0$  the graph $G(t+1)$ is constructed from $G(t)$ in the following way.
Add a new vertex $v_{t+3}$, then:
\begin{itemize}
    \item  with probability $1-\alpha$ connect $v_{t+3}$ to a vertex chosen by affine preferential attachment with attractiveness $\delta$: for all vertices $u\in V_t$, connect $v_{t+3}$ to $u$ with probability with a vertex $v_r\in V_{t}$,
    \begin{equation}
    \Prob(v_{t+3}\sim u)=(1-\alpha)\frac{d_t(ur)+\delta}{\sum_{i=1}^{t+2}(d_t(v_i)+\delta)}=(1-\alpha)\frac{d_t(u)+\delta}{2|E_t|+\delta(t+2)},
    \end{equation}
    where $d_t:=\{d_t(v_i),1\leq i\leq t+2\}$ is the degree vector in $G(t)$.
    \item With probability $\alpha$ connect $v_{t+3}$ with and edge in $E_t$, it chosen uniformly at random, that is, if $\{u,v\}\in E_t$
    \begin{equation}
    \Prob(v_{t+3}\sim u, v_{t+3}\sim v )=\frac{\alpha}{|E_t|}.
    \end{equation}
\end{itemize}

It is worth noticing that, when $\delta=0$, the recursion can also be defined this way:
\begin{itemize}
    \item The new vertex $v_{t+3}$ is connected to a single vertex $u$, chosen with probability proportional to its degree:
   \begin{equation}
    \Prob(v_{t+3}\sim u)=\frac{d_t(u)}{\sum_{i=1}^{t+2}d_t(v_i)}=\frac{d_t(u)}{2|E_t|}
    \end{equation}

    \item Given that at step 1 the vertex $v_r$ has been connected to $u$, with probability $1-\alpha$ do not do anything else, but with probability $\alpha$ 
    also connect it to a uniform neighbour of $u$.  
\end{itemize}

Note that both points of view on the triangle formation step are equivalent, since the probability of linking $v_{t+3}$ to both elements of an edge $\{u,v\}$ is equal to (separating depending on which of $u$ and $v$ is chosen first)
\[\left(\frac{d_t(u)}{2|E_t|}\right)\left(\frac{1}{d_t(u)}\right)+\left(\frac{d_t(v)}{2|E_t|}\right)\left(\frac{1}{d_t(v)}\right)=\frac{1}{|E_t|}.\]


This model reflect a reasonably natural social behaviour. When a new individual A arrives, it becomes friends of a popular individual B, and might also become friends with one of B's own friends.

\textbf{\textit{Notation.}} Given two sequences of real numbers $(x_t)_t$ and $(y_t)_t$ we write $x_t=O(y_t)$ provided that, for all large enough $t$, there is a positive constant $C$ such that $|x_t|\leq Cy_t$, while we write $x_t=\Theta(y_t)$ provided that, for all large enough $t$,  there are constants $C,c>0$ such that $cy_t\leq x_t\leq Cy_t$. Moreover, when we say that these bounds are uniform in another variable, say $n$, this means the constants $c$ and/or $C$ can be chosen to not depend of $n$.  

\section{Main Results}
\subsection{Power law degree distribution}
In this section we study the empirical degree distribution of vertices in our model. Let us write, for $t\geq 0$, $N(\ell;t)$ and $N(\ell,m;t)$ for respectively the number of vertices with degree $\ell$ in $G(t)$ and the number of edges with degrees $\ell$ and $m$ in $G(t)$, respectively, that is
\begin{align*}
N(\ell;t)=&\sum_{i=1}^{|V_t|}\mathbbm{1}_{\{d_t(v_i)=\ell\}}, \quad \text{ and}\\
N(\ell,m;t)=&\sum_{i=1}^{|V_t|}\sum_{j>1}^{|V_t|}\mathbbm{1}_{\{d_t(v_i)=\ell,d_t(v_j)=m,v_i\sim v_j\}}.
\end{align*}
Note that the total number of vertices $|V_t|=t+2$ for all $t\geq 0$.
\begin{thm}\label{Teo1}
Let
\begin{equation}\label{A}
    A=\frac{[(\alpha+1)^2+\delta\alpha]}{(\alpha+1)[2(\alpha+1)+\delta]},
\end{equation}
and for $\ell\geq 1$,
\begin{equation}\label{Aele}
    A_{\ell}=\frac{[(\alpha+1)^2+\delta\alpha]+\delta(1+\alpha)(1-\alpha)/\ell}{(\alpha+1)[2(\alpha+1)+\delta]}.
\end{equation}
For any $\epsilon>0$
\begin{equation}
    \lim_{t\rightarrow\infty}\Prob\Big(\Big|\frac{N(\ell,t)}{t}- p(\ell)\Big|<\epsilon\Big)=1,
\end{equation}
where
\begin{align}\label{pele}
    p(\ell)=
    \begin{cases}
    \frac{1-\alpha}{A_1+1}, & \text{if $\ell=1$}\\
    \frac{A_1+\alpha}{(A_1+1)(2A_2+1)}, & \text{if $\ell=2$}\\
    \frac{A_1+\alpha}{A_1A_{\ell}\prod_{j=1}^{\ell} \big(1+\frac{1}{jA_j}\big)}, & \text{if $\ell\geq 3$}.
    \end{cases}
\end{align}
Moreover, 
$$
p(\ell)=\Theta(\ell^{-\gamma}),\quad \text{ with } \gamma={1+\frac{1}{A}}.
$$
\end{thm}

\begin{obs}
Since $\gamma=1+\frac{1}{A}$, then there are three different regimes in this model, 
\begin{itemize}
    \item $\gamma>3$ when $A<1/2$, \item $\gamma=3$ when $A=1/2$ \item and $\gamma<3$ when $A>1/2$.
\end{itemize}
Note that if $\alpha=0$, that is, triangles are not included, then $A=1/(2+\delta)$. So $A$ can be smaller than, equal to, or larger than $1/2$ depending on the sign of $\delta$. In this case $p(\ell)=\Theta(\ell^{-3-\delta}),$ which corresponds to the model studied by M\'ori in \cite{Mori}.
The case $\delta=0$ and $\alpha=0$  corresponds to the Barab\'asi-Albert model with $m=1$,  numerically studied in \cite{BarabasiAlbert99} and analytically studied in \cite{Bollobas2001}. 
\end{obs}

\begin{obs}
The three different regimes can also be described as follows:
\begin{itemize}
    \item $A<1/2$ when $\delta>0$ and $\alpha<1$, 
    \item $A=1/2$ when $\delta=0$, or $\alpha=1$, 
    \item $A>1/2$ when $\delta<0$ and $\alpha<1$.
\end{itemize} 
\end{obs}

\subsection{Clustering coefficients}
Given a graph, a clustering coefficient is a measure of the degree to which its vertices tend to cluster. Two of most known such measures are the \textit{global} and the \textit{average local} clustering coefficients. The global clustering coefficient measures the proportion of triangles with respect to the total number of possible triangles, while the average local clustering coefficient measures the proportion, of the proportion of triangles between the neighbours of each vertex, with respect to the number of vertices.  
While these two objects are related, they can at times have very different behaviours. 
In this section we gives some asymptotic properties of both these clustering coefficients for our model.

\subsubsection{Average local clustering coefficient}

The local clustering coefficient at $v\in V_t$ is given by
\begin{equation}\label{Ci}
C_t(v)=
\begin{cases}
    \frac{|E_t(v)|}{\binom{d_t(v)}{2}}& \text{ if } d_t(v)\geq 2\\
    0& \text{ otherwise, }
\end{cases}
\end{equation}
where $|E_t(v)|$ is the number of edges between neighbours of $v$ in $G(t)$, or in other words the number of triangles in $G(t)$ which include $v$. The average local clustering $\mathcal{C}_2(t)$ of $G(t)$ is then defined by
\begin{equation}\label{localClustering}
\mathcal{C}_1(t)=\frac1{|V_t|}\sum_{i=1}^{|V_t|}C_t(v_i).
\end{equation}

\begin{thm}[The average local clustering coefficient does not vanish]\label{TeoLocalClustering}
    The average local clustering coefficient $\mathcal{C}_1(t)$ satisfies  
    \begin{equation}
    \lim_{t\rightarrow\infty}\Prob\Big(\mathcal{C}_1(t)\geq\alpha p(2)\Big)=1,
    \end{equation}
    with $p(2)$ defined in \eqref{pele}.
\end{thm}

\subsubsection{Global clustering coefficient}


Let $T_t$ be the number of triangles and $C_t$ the number of connected triples in $G(t)$, where a triangle in $G(t)$ is a set of three vertices ${u,v,w}$ in $V_t$ that are mutually connected by edges, and  
  
 \begin{equation}\label{Pt}
 C_t=\sum_{i=1}^{|V_t|}\binom{d_t(v_i)}{2}=\sum_{\ell=2}^{t+1}N(\ell;t)\frac{\ell(\ell-1)}{2}.
\end{equation}
The global clustering coefficient $\mathcal{C}_2(t)$ of $G(t)$ is defined by
\begin{equation}\label{globalClustering}
    \mathcal{C}_2(t)=\frac{3T_t}{C_t}.
\end{equation}

\begin{obs}
It is not difficult to see that an equivalent formulation for $\mathcal{C}_2(t)$ is given by
\[
\mathcal{C}_2(t)=\frac{\sum_{i=1}^{|V_t|} \binom{d_t(v_i)}{2} C_t(v_i)}{\sum_{i=1}^{|V_t|} \binom{d_t(v_i)}{2}},
\]
with $C_t(v_i)$ as in (\ref{Ci}). This lets us see that $\mathcal{C}_2(t)$ is a weighted average of the $C_t(v_i)$ terms, which gives higher weight than $\mathcal{C}_1(t)$ to vertices with higher degree.
\end{obs}

To study the expectation of global clustering coefficient we will need to study the expectation of the number of connected triples, $C_t$. The authors of \cite{Ostroumova}
say that the value $C_t$ in scale-free graphs is determined by the exponent $\gamma$ of the power law, and  write the following: ``$C_t=\sum_{\ell=2}^{t+1}N(\ell;t)\frac{\ell(\ell-1)}{2}\propto \sum_{\ell=2}^{t+1} t \ell^{2-\gamma}$. Therefore if $\gamma>3$, then $C_t$ is linear in $t$. However, if $\gamma\leq 3$, $C_t$ is superlinear''.  
\begin{obs}For the model  described in Section \ref{TheModel}, $\gamma=1+\frac{1}{A}$. Lemma \ref{lemEPt} in the proof's section confirms the findings in \cite{Ostroumova}, but just in terms of the expected value. 
\end{obs}

Our result for the expected value of the global clustering coefficient is the following.

\begin{thm}\label{ThmExpC2}
   Consider the model described in Section \ref{TheModel} and $A$ as in (\ref{A}). Then, the expected value of the global clustering coefficient is given by 
   \begin{equation}
    \E[\mathcal{C}_2(t)]=
     \begin{cases}
    \Theta(1), & \text{if $A<\frac12$}\\
    \Theta(ln(t)^{-1}), & \text{if $A=\frac12$}\\
    \Theta(t^{1-2A}), & \text{if $A>\frac12$}
    \end{cases}
\end{equation}
\end{thm}

\begin{obs}  In \cite{Oliveira}, the authors study $C_t$ for the model of \cite{HolmeKim02}. They used martingale theory and a combination of  Azuma-Hoeffding and Freedman's concentration inequalities, together with a bootstrap argument to get a weak-law of large numbers for $C_t$. This method unfortunately isn't sufficient for our model. The reason is that these two inequalities need that the jumps of the martingales in question are bounded, but that is not satisfied when the number of edges is random which is the case for our model. Thus, other techniques are needed to obtain concentration results for $C_t$ around its mean. 
\end{obs}

\section{Proofs}\label{RecExp}
\subsection{Degree distribution}
To prove Theorem \ref{Teo1} we start with Lemma \ref{Lem1} which provides recursive formulae for the expected value of the number of vertices with a given degree in $G(t+1)$. Lemma \ref{Lemap(l)} then uses this to establish than the expected proportion of vertices with a given degree has the wanted limit as $t\rightarrow\infty$ of the expected, and finally, Theorem \ref{Teo1} is proved at the end by showing that these proportions are close to their expectation.

\begin{lem}\label{Lem1}
The following estimates hold:

 \begin{align*}   
    \E[N(1;t+1)]&=(1-\alpha)+\left(1-\frac{A_1}{t}\right)\E[N(1;t)]+O\left(\frac{1}{\sqrt{t}}\right),\\
    \E[N(2;t+1)]&=\alpha+\frac{A_1\E[N(\ell-1;t)]}{t}+\left(1-\frac{2A_2}{t}\right)\E[N(2;t)]+O\left(\frac{1}{\sqrt{t}}\right),\\
    \E[N(\ell;t+1)]&=(\ell-1)\frac{A_{\ell-1}\E[N(\ell-1;t)]}{t}+\left(1-\frac{\ell A_{\ell}}{t}\right)\E[N(\ell;t)]+O\left(\frac{1}{\sqrt{t}}\right), & \text{for $\ell>2.$}
\end{align*}
where $A_{\ell}$ is given by (\ref{Aele}), and the error term  $O\left(\frac{1}{\sqrt{t}}\right)$ is uniform in $\ell$.
\end{lem}

\begin{proof}
Let, for $s\geq0,$ $\mathcal{G}_s$ be the $\sigma$-field generated by the appearance of edges up to time $s$. We study the distribution of the increment $N(\ell;t+1)-N(\ell;t)$, conditionally of $\mathcal{G}_t$, separating the cases $\ell=1$, $\ell=2$ and $\ell> 2$.

Let us start by considering the case $\ell>2$. Note that $N(\ell;t+1)-N(\ell;t)$ can be any integer from $-2$ to $+2$.

The increment is equal to 2 (respectively $-2$)  when the added vertex at time $t+1$, $v_{t+3}$, is connected with and edge such that, its both incident vertices  have degrees $\ell-1$ (respectively $\ell$). This happens with probability $\alpha/|E_t|$ for each such edge. 

The increment is equal to $1$ (respectively $-1$) when one of the following events occurs. The first one is when $v_{t+3}$ is connected with an edge such that, just one of its incident  vertices  has degree $\ell-1$ (respectively $\ell$). That happens with probability $\alpha/|E_t|$. The second one is when $v_{t+3}$  is connected with an existing vertex with degree $\ell-1$ (respectively $\ell$). This happens with probability $(1-\alpha)(\ell-1+\delta)/(2|E_t|+\delta(t+2))$ (respectively $(1-\alpha)(\ell+\delta)/(2|E_t|+\delta(t+2))$) per such vertex.

Finally, the increment is $0$ in all other situations. Using all of this we achieve 
\begin{align}\label{Nl3}
N(\ell;t+1)=
    \begin{cases}
    N(\ell;t)+2,  & w.p\quad p_2 \\
    N(\ell;t)-2,  & w.p\quad q_2 \\
    N(\ell;t)+1,  & w.p\quad p_1 \\
    N(\ell;t)-1,  & w.p\quad q_1 \\
    N(\ell;t),  & w.p\quad 1-p_1-p_2-q_1-q_2 \\
    \end{cases}
\end{align}
where $w.p$ means "with probability" and
\begin{align}
   p_1&=\frac{\alpha}{|E_t|}\sum_{k=2,k\neq \ell-1,\ell}^{t+2}N(\ell-1,k;t)+\frac{(1-\alpha)(\ell-1+\delta)N(\ell-1;t)}{2|E_t|+\delta(t+2)}\\
   q_1&=\frac{\alpha}{|E_t|}\sum_{k=2,k\neq \ell-1,\ell}^{t+2}N(\ell,k;t)+\frac{(1-\alpha)(\ell+\delta) N(\ell;t)}{2|E_t|+\delta(t+2)}\\
   p_2&=\frac{\alpha}{|E_t|}N(\ell-1,\ell-1;t)\\
   q_2&=\frac{\alpha}{|E_t|}N(\ell,\ell;t).
\end{align}

Taking the conditional expectations of (\ref{Nl3}) given $\mathcal{G}_t$, and noticing that $N(\ell-1,\ell;t)=N(\ell,\ell-1;t)$ (since edges are undirected), we get
\begin{align}\label{ENcaso3}
 \E[N(\ell;t+1)]&-\E[N(\ell;t)]=\E\{\E[N(\ell;t+1)-\E[N(\ell;t)]|\mathcal{G}_t]\}\\
    &= \E\left\{\frac{\alpha}{|E_t|}\left[2N(\ell-1,\ell-1;t)+\sum_{k=2,k\neq\ell-1}^{t+2}N(\ell-1,k;t)\right]\right\}\nonumber\\
    &- \E\left\{\frac{\alpha}{|E_t|}\left[2N(\ell,\ell;t)+\sum_{k=2,k\neq\ell}^{t+2}N(\ell,k;t)\right]\right\}\nonumber\\
    &+  \E\left\{\frac{(1-\alpha)(\ell-1+\delta)N(\ell-1;t)}{2|E_t|+\delta(t+2)}-\frac{(1-\alpha)(\ell+\delta) N(\ell;t)}{2|E_t|+\delta(t+2)}\right\}.
\end{align}
Now observe the following:
\begin{align}\label{l-1caso3}
& 2N(\ell-1,\ell-1;t)+\sum_{k=2,k\neq\ell-1}^{t+2}N(\ell-1,k;t)\nonumber\\
&=N(\ell-1,\ell-1;t)+\sum_{k=2}^{t+2}N(\ell-1,k;t)=(\ell-1)N(\ell-1;t).
\end{align}
The reason for this is that $(\ell-1)N(\ell-1;t)$ counts the number of edges connected to at least one vertex with degree $\ell-1,$ and we split it based on the degree of the other vertex - noticing moreover that if both vertices have degree $\ell$, then the edge is counted twice. 

Similarly we have 
\begin{align}\label{lcaso3}
    2N(\ell,\ell;t)+\sum_{k=2,k\neq\ell}^{t+2}N(\ell,k;t)=\ell N(\ell;t).
\end{align}
Substituting \eqref{l-1caso3} and \eqref{lcaso3} in \eqref{ENcaso3}, a few calculations then lead to
\begin{align}\label{ENl3}
    \E[N(\ell;t+1)]=&  \E[N(\ell;t)]+\E\left[\frac{\alpha(\ell-1)N(\ell-1;t)}{|E_t|} + \frac{(1-\alpha)(\ell-1+\delta)N(\ell-1;t)}{2|E_t|+\delta(t+2)}\right]\nonumber\\
     &-\E\left[\frac{\alpha\ell N(\ell;t)}{|E_t|} + \frac{(1-\alpha)(\ell+\delta) N(\ell;t)}{2|E_t|+\delta(t+2)}\right]\nonumber.\\
\end{align}

Similarly, when $\ell=2$, conditioning on $\mathcal{G}_t$, the probability distribution of $N(\ell; t + 1)$ is given by 
\begin{align}\label{Nl2}
N(\ell;t+1)=
    \begin{cases}
    N(\ell;t)+2,  & w.p\quad \hat{p}_2 \\
    N(\ell;t)+1,  & w.p\quad \hat{p}_1 \\
    N(\ell;t)-1,  & w.p\quad \hat{q}_1 \\
    N(\ell;t),  & w.p\quad 1-\hat{p}_1-\hat{p}_2-\hat{q}_1 \\
    \end{cases}
\end{align}
where
\begin{align}
   \hat{p}_1&=\alpha\left(1-\frac{1}{|E_t|}\sum_{k=2}^{t+2}N(2,k;t)-\frac{1}{|E_t|}\sum_{k=3}^{t+2}N(1,k;t)\right)+\frac{(1-\alpha)(1+\delta)N(1;t)}{2|E_t|+\delta(t+2)}\\
   \hat{q}_1&=\frac{\alpha}{|E_t|}N(2,2;t)+\frac{(1-\alpha)(2+\delta)}{2|E_t|+\delta(t+2)}N(2;t)\\
   \hat{p}_2&=\frac{\alpha}{|E_t|}\sum_{k=3}^{t+2} N(1,k;t).
\end{align}
Taking the expectation of the conditional expectations found in \eqref{Nl2}, we get
\begin{align}\label{ENl2}
    \E[N(2;t+1)]=&  \alpha+\E[N(2;t)]+\E\left[\frac{\alpha N(1;t)}{|E_t|} + \frac{(1-\alpha)(1+\delta)N(1;t)}{2|E_t|+\delta(t+2)}\right]\nonumber\\
     &-\E\left[\frac{2\alpha N(2;t)}{|E_t|} + \frac{(1-\alpha)(2+\delta) N(2;t)}{2|E_t|+\delta(t+2)}\right]\nonumber\\
\end{align}
Finally,  for the case $\ell=1$
\begin{align}\label{Nl1}
N(\ell;t+1)=
    \begin{cases}
    N(\ell;t)+1,  & w.p\quad \tilde{p}_1 \\
    N(\ell;t)-1,  & w.p\quad \tilde{q}_1 \\
    N(\ell;t),  & w.p\quad 1-\tilde{p}_1-\tilde{q}_1 \\
    \end{cases}
\end{align}
where
\begin{align}
   \tilde{p}_1&=(1-\alpha)\left(1-\frac{(1+\delta)N(1;t)}{2|E_t|+\delta(t+2)}\right)\\
   \tilde{q}_1&=\frac{\alpha}{|E_t|}\sum_{k=2}^{t+2}N(1,k;t),
   \end{align}
and taking the expectation of the conditional expectations in  \eqref{Nl1} given $\mathcal{G}_t$,  and after few calculations we obtain 
\begin{align}\label{ENl1}
    \E[N(1;t+1)]=&  (1-\alpha)+\E[N(1;t)]
     -\E\left[\frac{\alpha N(1;t)}{|E_t|} + \frac{(1-\alpha)(1+\delta) N(1;t)}{2|E_t|+\delta(t+2)}\right]\nonumber\\
\end{align}
To finish the proof, we need to show the following:
\begin{equation}\label{27}
\E\left[\frac{N(\ell,t)}{|E_t|}\right] = \frac{\E[N(\ell,t)]}{t(\alpha+1)}+O\left(\frac{1}{\sqrt{t}}\right)
\end{equation}
and
\begin{equation}\label{28}
\E\left[\frac{N(\ell,t)}{2|E_t|+\delta(t+2)}\right] = \frac{\E[N(\ell,t)]}{t(2\alpha+2)+\delta t}+O\left(\frac{1}{\sqrt{t}}\right)
\end{equation}
where the $O(\frac1{\sqrt{t}})$ terms are uniform in $\ell.$ These are then easily subbed in (\ref{ENl1}), (\ref{ENl2}) and (\ref{ENl3})  to yield the Lemma. Since the proofs of (\ref{27}) and (\ref{28}) are essentially the same, we focus on (\ref{27}). Recall that $|E_t|$, the number of edges in $G(t)$ is a random variable. By construction, we can write it for all $t\geq 0$ as
\[
|E_t|=1+\sum_{i=1}^t X_i,
\]
where $(X_k,k\geq 1)$ is an i.i.d. sequence of random variables equal to $1$ with probability $1-\alpha$ and to $2$ with probability $\alpha$. In particular we have 
$\E[|E_t|]=1+t(\alpha+1)$ and $\mathrm{Var}[|E_t|]=t\alpha(1-\alpha).$

Now, we have
\[
    \E\left[\frac{N(\ell,t)}{|E_t|}\right] = \E\left[\frac{N(\ell,t)}{t(\alpha+1)}\right]+ \E\left[N(\ell,t)\left(\frac1{|E_t|}-\frac1{t(\alpha+1)}\right)\right],
\]
and we deduce from $N(\ell,t)\leq 2+t$ that
\[\left|\E\left[\frac{N(\ell,t)}{|E_t|}\right] - \E\left[\frac{N(\ell,t)}{t(\alpha+1)}\right]\right| \leq (2+t)\E\left[\frac{||E_t|-t(\alpha+1)|}{t(\alpha+1)|E_t|}\right],\]
and since $|E_t|\geq t+1$, we deduce
\[\left|\E\left[\frac{N(\ell,t)}{|E_t|}\right] - \E\left[\frac{N(\ell,t)}{t(\alpha+1)}\right]\right| \leq \frac{2}{t(\alpha+1)}\E[||E_t|-t(\alpha+1)|].\]


Finally,  using the Cauchy-Schwarz inequality, we have 
\begin{align*}
    \E[||E_t|-t(\alpha+1)|]&\leq \sqrt{\E\left[\left(|E_t|-t(\alpha+1)\right)^2\right]} \\
    &\leq \sqrt{\E\left[(|E_t|-t(\alpha+1))^2\right]+\mathrm{Var}[|E_t|]} \\
    &= \sqrt{1+t\alpha(1-\alpha)},
\end{align*}
and \eqref{27} follows readily.
\end{proof}

The following result states the limit behaviour of the expect proportion of vertices with a given degree in our model.  

\begin{lem}\label{Lemap(l)}
We have 
\begin{equation}
    \frac{\E[N(\ell;t)]}{t}\underset{t\rightarrow\infty}\longrightarrow p(\ell),
\end{equation}
where $p(\ell)$ was defined in \eqref{pele}.
Moreover,
\begin{align*}
    p(\ell)=\Theta(\ell^{-1-\frac{1}{A}}),
\end{align*}  where $A$ and $A_{\ell}$, $\ell\geq 1$, are  given by (\ref{A}) and (\ref{Aele}), respectively. 
\end{lem}

\begin{proof}
We start with the case $\ell=1$. Iterating one time the corresponding equation in Lemma \ref{Lem1}, we get that  $\E[N(1;t+1)]$ is equal to
\[
 \left(1-\frac{A_1}{t}\right)\left(1-\frac{A_1}{t-1}\right)\E[N(1;t-1)]+\left(1-\frac{A_1}{t}\right)(1-\alpha)+(1-\alpha)
+O\left(\sum_{i=t-1}^ t \sqrt{\frac{1}{i}}\right).    
\]

Iterating $t-1$ times, 
we then end up with

\[
\E[N(1;t+1)]=
\E[N(1;1)]\prod_{r=1}^{t} \left(1-\frac{A_1}{r}\right)+(1-\alpha)\sum_{s=1}^t\prod_{r=s+1}^t\left(1-\frac{A_1}{r}\right)+O\left(\sum_{i=1}^ t \sqrt{\frac{1}{i}}\right),
\]

Observe that, since $0<A_1<1$ for all choices of $\delta>-1$, we have $0<\prod_{r=1}^{t} \left(1-\frac{A_1}{r}\right)<1$, so the first time of the above sum is bounded and is a $o(t).$ The final term is, by comparison with an integral, $O(\sqrt{t}),$  hence to find the limit of $\frac{\E[N(1;t+1)]}{t}$ we only need to investigate the middle term $(1-\alpha)\sum_{s=1}^t\prod_{r=s+1}^t\left(1-\frac{A_1}{r}\right).$

Write $\prod_{r=s+1}^t\left(1-\frac{b}{r}\right)=exp\left(\sum_{r=s+1}^t\ln\left(1-\frac{b}{r}\right)\right)$. If $b\in(0,2)$ is fixed then, by Taylor expansion, for $r\geq 2$ $\ln\left(1-b/r\right)=-b/r-z_r$, with $z_r=O(b/r)^2$. Using $e^{-x}>1-x$ for any $x\in \R$, and 
\[
\int_{m}^{N+1}\left(\frac{1}{x}\right)^j dx\leq\sum_{i=m}^{N}\left(\frac{1}{i}\right)^j \leq\left(\frac{1}{m}\right)^{j}+ \int_{m}^{N}\left(\frac{1}{x}\right)^j dx,
\]
for any $m\geq 1$, $j\geq 1$ and $N>m$, then we obtain that for $s,t\in \N,$
\begin{equation}\label{producto}
\prod_{r=s+1}^t\left(1-\frac{b}{r}\right)=\left(\frac{s}{t}\right)^b\left(1+O\left(\frac{1}{s}\right)\right),
\end{equation}
where the $O$ term is uniform in $s$ and $t$.
In consequence, 
\[
\sum_{s=1}^t\prod_{r=s+1}^t\left(1-\frac{A_1}{r}\right)=\sum_{s=1}^t \Big(\frac{s}{t}\Big)^{A_1}\Big(1+O\Big(\frac{1}{s}\Big)\Big).
\]

where 
$$
\sum_{s=1}^t \Big(\frac{s}{t}\Big)^{A_1}\sim \int_0^t \Big(\frac{s}{t}\Big)^{A_1}ds=\frac{t}{A_1+1},
$$
and
$$
\sum_{s=1}^t \Big(\frac{s}{t}\Big)^{A_1}O\Big(\frac{1}{s}\Big)=O(1),
$$
as $t\rightarrow\infty$.
Therefore

\begin{equation}\label{p1}
\frac{\E[N(1,t)]}{t}\underset{{t\rightarrow \infty}}\longrightarrow \frac{(1-\alpha)}{A_1+1}=p(1).
\end{equation}
We use induction to treat the $\ell\geq 2$ cases. First, these cases of Lemma \ref{Lem1} can be grouped as
\begin{equation}\label{Eele}
\E[N(\ell,t+1)]=\left(1-\frac{\ell A_{\ell}}{t}\right)\E[N(\ell,t)]+g(\ell,t),
\end{equation}

where
\begin{align}\label{gt}
    g(\ell,t)=
    \frac{(\ell-1)A_{\ell-1}\E[N(\ell-1,t)]}{t}+\alpha\mathbbm{1}_{\ell=2}+O\left(\sqrt{\frac1{t}}\right)
\end{align}
By induction, assuming that $\E[N(\ell-1,t)]/t$ converges to $p(\ell-1),$ then $g(\ell,t)$ has limit $G(\ell)$ as $t\rightarrow\infty$, defined by
\begin{align}\label{g}
    G(\ell)=(\ell-1)A_{\ell-1}p(\ell-1) + \alpha\mathbbm{1}_{\ell=2}.
\end{align}

Now we follow the ideas of the proof of Theorem 4.1.2 in \cite{durrett_2006}. 
Iterating  (\ref{Eele}) $t-1$ times when $\ell=2$, and $t-(\ell-2)$ times when $\ell>2$,  we get
\begin{align*}
    &\E[N(\ell;t+1)]=\\
    &\begin{cases}
    \E[N(2,1)]\prod_{r=1}^t\left(1-\frac{2 A_{2}}{r}\right)+\sum_{s=1}^{t}g(2,s)\prod_{r=s+1}^t\left(1-\frac{2 A_{2}}{r}\right) , & \text{if } \ell=2\\
    \sum_{s=\ell-2}^{t}g(\ell,s)\prod_{r=s+1}^t\left(1-\frac{\ell A_{\ell}}{r}\right) , & \text{if } \ell>2,
    \end{cases}
\end{align*}
with $\prod_{r=t+1}^t\left(1-\frac{\ell A_{\ell}}{r}\right)=1$.

Note that $A_{\ell}<1$ for $\delta>-1$,  $\E[N(2;1)]=3\alpha+(1-\alpha)$ and for $\ell>2$, $\E[N(\ell;t)]>0$ only for  $t\geq\ell-1$.
Moreover,  $0<\prod_{r=1}^{t} \left(1-\frac{2A_2}{r}\right)<1$, so $\frac{1}{t}\E[N(2,1)]\prod_{r=1}^{t} \left(1-\frac{2A_2}{r}\right)\rightarrow 0$ as $t$ goes to infinity.

The main part of the induction step will be to combine \eqref{producto} and the limit of $g(\ell,t)$ and obtain
\begin{equation}\label{beforelimitE}\sum_{s=a}^{t}g(\ell,s)\prod_{r=s+1}^t\left(1-\frac{\ell A_{\ell}}{r}\right)\underset{t\to\infty}\sim \sum_{s=a}^{t}G(\ell)\left(\frac{s}{t}\right)^{\ell A_{\ell}}\underset{t\to\infty}\sim \frac{tG(\ell)}{\ell A_\ell +1},
\end{equation}
where $a=1$ when $\ell=2$ and $a=\ell-2$ for $\ell>2$, giving

\begin{equation}\label{limitE}
\frac{\E[N(\ell,t)]}{t}\underset{t\to\infty}\longrightarrow \frac{G(\ell)}{\ell A_{\ell}+1}.
\end{equation}

First, write
\[P_t=\sum_{s=a}^{t}|g(\ell,s)-G(\ell)|\prod_{r=s+1}^t\left(1-\frac{\ell A_{\ell}}{r}\right)\]
and
\[Q_t=\sum_{s=a}^{t}G(\ell)\left|\prod_{r=s+1}^t\left(1-\frac{\ell A_{\ell}}{r}\right)-\left(\frac{s}{t}\right)^{\ell A_{\ell}}\right|,\]
and let us show that both $P_t$ and $Q_t$ are $o(t)$ as $t\to\infty,$ which is enough to deduce \eqref{beforelimitE}.

For $P_t$, write $P_t\leq \sum_{s=a}^{t}|g(\ell,s)-G(\ell)|.$ Since $|g(\ell,s)-G(\ell)|=o(1)$, summation of asymptotics gives that $P_t=o(t).$ For $Q_t,$  we have
\[Q_t=G(\ell)\sum_{s=a}^t \left(\frac{s}{t}\right)^{\ell A_{\ell}}O\left(\frac1{s}\right)=\frac{1}{t^{\ell A_{\ell}}}\sum_{s=a}^t O(s^{\ell A_{\ell}-1}).\]
Summation of asymptotics then gives us $Q_t=\frac{1}{t^{\ell A_{\ell}}}O(t^{\ell A_{\ell}})=O(1)=o(t),$ completing the proofs of \eqref{beforelimitE} and \eqref{limitE}.

To finish the induction, we need to check that $\frac{G(\ell)}{\ell A_{\ell}+1}$ is equal to $p(\ell)$ as defined in the statement. The $\ell=2$ case is immediate since by (\ref{p1}) and (\ref{g}), $p(2)=\frac{A_1p(1)+\alpha}{2A_2+1}=\frac{A_1+\alpha}{(A_1+1)(2A_2+1)}.$ For the further cases, notice first from the definition of $p(\ell)$ \eqref{pele} that
\begin{equation}\label{pl}
p(\ell)=p(\ell-1)\frac{(\ell-1)A_{\ell-1}}{\ell A_\ell}\frac1{1+\frac1{\ell A_{\ell}}},
\end{equation}
then by definition of $G(\ell)$ in (\ref{g}), we have
\[\frac{G(\ell)}{\ell A_{\ell}+1}=\frac{(\ell-1)A_{\ell-1}p(\ell-1)}{\ell A_{\ell}+1}=p(\ell),\]
and this ends the induction.
Finally, we check the asymptotic for $p(\ell).$ Note that $jA_j=jA+B$, where 
\begin{equation}\label{AB}
    A=\frac{(\alpha+1)^2+\delta\alpha}{(\alpha+1)[2(\alpha+1)+\delta]} \quad\text{ and }\quad B=\frac{\delta(1+\alpha)(1-\alpha)}{(\alpha+1)[2(\alpha+1)+\delta]},
\end{equation}
and therefore $A_{\ell}\sim A$. Moreover,

\begin{align}\label{expterm}
\prod_{j=1}^{\ell} \Big(1+\frac{1}{jA_{j}}\Big)&=\exp^{\sum_{j=1}^{\ell}\ln \Big(1+\frac{1}{jA_{j}}\Big)}\nonumber\\
&=\Theta \left(\exp^{\int_{1}^{\ell}\ln \Big(1+\frac{1}{xA+B}\Big) dx} \right) \nonumber \\
&=\Theta\big(\ell^{\frac{1}{A}}\big).    
\end{align}




Therefore we conclude
$$
p(\ell)=\Theta\big(\ell^{-1-\frac{1}{A}}\big)
$$

\end{proof}

The proof of Theorem \ref{Teo1} uses the Azuma-Hoeffding inequality, which we state here:

\begin{thm}[Inequality
of Azuma \cite{Azuma} and Hoeffding \cite{hoeffding1963probability}]\label{Azuma}
 Let $(X_s,s\in \{0,1,\ldots,t\})$ be a martingale with respect to certain filtration $(\mathcal{F}_s,s\in \{0,1,\ldots,t\})$ and which satisfies $|X_s - X_{s-1}|\leq c$ a.s. for $1 \leq s \leq t$. Then
\[\Prob(|X_t - X_0| \geq x) \leq \exp(-x^2/2c^2 t).
\]
\end{thm}

\noindent\textit{Proof of Theorem \ref{Teo1}.} Define $W_s = \E [N(\ell, t) \mid \mathcal{G}_s]$ for $s\leq t$, where $\mathcal{G}_s$ is the $\sigma$-field generated by the appearance of edges up to time $s$ in the random graph model. Then $W_s$ is a martingale with respect to $\mathcal{G}_s$ and
$|W_s-W_{s-1}|\leq 5$. This last observation happens because in one step a new vertex is added and it is attached to at most two existing vertices, which implies that  the degree of at most four existing vertices at time $s-1$ and the vertex added at time $s$, will be affected. Therefore, applying the Azuma-Hoeffding inequality with $W_t = N(\ell, t)$,
$W_0 = \E[N(\ell, t)]$ and taking $x =\sqrt{t \ln(t)}$, we obtain
\begin{equation}\label{TeoDegree}
    \Prob\left(\left|\frac{N(\ell, t)}{t}-\frac{ \E[N(\ell, t)]}{t}\right|\geq\sqrt{\frac{\ln t}{t}}\right)\leq t^{-1/50}.
\end{equation}

Thus by Lemma \ref{Lemap(l)}, as $t\rightarrow \infty$ we have $N(\ell, t)/t \rightarrow p(\ell)$ in probability.
\qed

\subsection{Average local clustering coefficient: proof of Theorem \ref{TeoLocalClustering}}
We aim to show that the probability that the average local clustering coefficient is bounded away from zero goes to one, as $t$ goes to infinity. We will follow some ideas from \cite{Oliveira}.

By construction of the model, at each time step a vertex is added and with probability $\alpha$ it is connected to an existing edge, creating in this way a triangle in $G(t)$. Let $T_{v_i}$ denote the number of triangles added at time $i-2$ for $i>2$ (the moment $v_i$ is added), then $T_{v_i}$ follows a Bernoulli distribution with $\Prob(T_{v_i}=1)=\alpha$. By Lemma \ref{Lemap(l)} we know that
$\E[N(\ell,t)]/t$ converges to $p(\ell)$ as $t$ goes to infinity.
Using this and by (\ref{TeoDegree}),  for all $t$ large enough we have 
\[
\Prob\Big(\frac{N(\ell,t)}{t}\geq \frac{p(\ell)}{2}-\sqrt{\frac{\ln t}{t}}\Big)\geq \Prob\Big(\frac{N(\ell,t)}{t}\geq \frac{\E[N(\ell,t)]}{t}-\sqrt{\frac{\ln t}{t}}\Big)\geq 1-t^{-1/50}.
\]
For $\ell\geq 2$, let $N^{\alpha}(\ell,t)$ be the number of vertices in $G(t)$ with degree $\ell$ which appeared with a triangle formation step: 
\[
N^{\alpha}(\ell,t)=\sum_{i=1}^{|V_t|} \mathbbm{1}_{\{d_t(v_i)=\ell, T_{v_i}=1\}}, \quad \ell\geq 2.
\]
Let $Bin(n,p)$ denote a generic binomial random variable with parameters $n$ and $p$. For any $x\in \{0,1,\dots,t+1\}$ we have

\begin{align*}
    \Prob(N^{\alpha}&(\ell,t)\leq x) \\
    &\leq \sum_{j=\frac{p(\ell)t}{2}-\sqrt{t\ln t}}^{t+2}\Prob(N^{\alpha}(\ell,t)\leq x\mid N(\ell,t)=j)\Prob(N(\ell,t)=j)+\Prob\Big(N(\ell,t)<\frac{p(\ell)t}{2}-\sqrt{t\ln t}\Big)\\
    &=\sum_{j=\frac{p(\ell)t}{2}-\sqrt{t\ln t}}^{t+2} \Prob(Bin(j,\alpha)\leq x)\Prob(N(\ell,t)=j)+\Prob\Big(N(\ell,t)<\frac{p(\ell)t}{2}-\sqrt{t\ln t}\Big)\\
    &\leq\;\Prob\Big(Bin\Big(\frac{p(\ell)t}{2}-\sqrt{t\ln t},\alpha\Big)\leq x\Big)+t^{-1/50}.
\end{align*}
By Chernoff's inequality applied to a Binomial random variable, see for example Theorem 2.1 in \cite{Janson}, we have that, for any $\ell\geq 2,$

\begin{equation}\label{Nalpha}
\Prob\left(N^{\alpha}(\ell,t)\leq \Big(\frac{p(\ell)t}{2}-\sqrt{t\ln t}\Big)\alpha-\sqrt{\alpha t\ln t}\right)
\leq \mathrm{exp}\left(\frac{-\alpha t\ln t}{2\alpha \left(\frac{p(\ell)t}{2}-\sqrt{t\ln t}\right)}\right)+t^{-1/50}\nonumber=o(1).
\end{equation}

Now, observe that for each vertex $v_i\in V_t$ such that, $T_{v_i}=1$, the number of edges between the neighbours of $v_i$ in $G(t)$, that is, $|E_t(v_i)|$, is at least 2. Using this when $d_t(v_i)\geq 2$, we have that the local clustering coefficient at $v_i$ given by (\ref{Ci}) satisfies
\[
C_t(v_i)\geq \frac{4}{d_t(v_i)(d_t(v_i)-1)},
\]
and by (\ref{localClustering})
\begin{align*}
\mathcal{C}_1(t)=\sum_{i=1}^{t+2}\mathbbm{1}_{\{T_{v_i}=1\}}\frac{C_t(v_i)}{t+2}&\geq  \sum_{i=1}^{t+2}  \mathbbm{1}_{\{T_{v_i}=1\}}\frac{4}{(t+2)d_t(v_i)(d_t(v_i)-1)}\\
&= \frac{4}{t+2}\sum_{\ell=2}^{t}\frac{N^{\alpha}(\ell,t)}{\ell(\ell-1)}\\
&\geq \frac{2N^{\alpha}(2,t)}{(t+2)}.
\end{align*}
Therefore, using this and by (\ref{Nalpha})
\begin{align*}
 1&=\Prob\Big(\mathcal{C}_1(t)\geq\frac{2N^{\alpha}(2,t)}{(t+2)}\Big)\\
    &\leq \Prob\left(\mathcal{C}_1(t)\geq\frac{2N^{\alpha}(2,t)}{(t+2)}, N^{\alpha}(2,t)\geq \Big(\frac{p(2)t}{2}-\sqrt{t\ln t}\Big)\alpha-\sqrt{\alpha t\ln t}\right) +o(1)\\
    &\leq \Prob\left(\mathcal{C}_1(t)\geq\frac{2\alpha\Big(\frac{p(2)t}{2}-2\sqrt{t\ln t}\Big)}{(t+2)}\right) +o(1),
\end{align*}
from which we have
\[
\underset{t\to\infty} \lim\Prob\left(\mathcal{C}_1(t)>\frac{2\Big(\frac{p(2)t}{2}-2\sqrt{t\ln t}\Big)\alpha}{(t+2)}\right)= 1.
\]
Since 
\[
\lim_{t\to \infty} \frac{2\Big(\frac{p(2)t}{2}-2\sqrt{t\ln t}\Big)\alpha}{(t+2)}=\alpha p(2),
\]
we deduce
\[
\lim_{t\rightarrow\infty} \Prob\Big(\mathcal{C}_1(t)\geq \alpha p(2)-\varepsilon\Big)=1
\]
for all $\varepsilon>0$, and hence
\[
\lim_{t\rightarrow\infty} \Prob\Big(\mathcal{C}_1(t)\geq \alpha p(2\Big)=1.
\]\qed

\subsection{Global clustering coefficient}

Recall that the global clustering coefficient $\mathcal{C}_2(t)$ of $G(t)$ can be written as
\[\mathcal{C}_2(t)=\frac{3T_t}{C_t},\]
where $C_t=\sum_{\ell=2}^{t+1}N(\ell;t)\frac{\ell(\ell-1)}{2}$
is the number of connected triples of $G(t).$
The behaviour of the numerator is easily studied: since it has a probability $\alpha$ of growing independently at each step, the law of large numbers gives us the following:

\begin{lem}\label{lemTt}We have $\E[T_t]=\alpha t$ for all $t\in\N$, and
\begin{equation}
    \frac{T_t}{t}\rightarrow \alpha
\end{equation}
almost surely as $t\to\infty.$
\end{lem}



Now we study the expected value of $C_t$.

\begin{lem}\label{lemEPt}
Recall that $A$ is given by \eqref{A}. We have
\begin{equation}
    \E[C_t]=
     \begin{cases}
    \Theta(t), & \text{if $2A<1$}\\
    \Theta(t\ln(t)), & \text{if $2A=1$}\\
    \Theta(t^{2A}), & \text{if $2A>1$}
    \end{cases}
\end{equation}
\end{lem}

\begin{proof}

Let $\tilde{C}_t:=\sum_{i=1}^{|V_t|}d_t^2(v_i)$ denote the sum of the squares of the degrees in $G(t)$. As in \cite{Oliveira}, we start by studying $\tilde{C}_t$  since $C_t=\frac{1}{2}\tilde{C}_t-|E_t|$.

Note that we can write
\begin{equation}\label{1}
    d_{t+1}(v_i)=d_t(v_i)+\Delta d_{t+1}(v_i),
\end{equation}
where $\Delta d_{t}(v_i):=d_{t}(v_i)-d_{t-1}(v_i)$ is the increment of $d_t(v_i)$ in one time step. For any $t\geq 1$,  $\Delta d_{t}(v_i)$  is a random variable taking values in $\{0,1\}$ and, for any fixed vertex $v\in V_t$
\begin{equation}\label {2}
    \Prob(\Delta d_{t+1}(v)=1\mid \mathcal{G}_t)=(1-\alpha)\frac{d_t(v)+\delta}{2|E_t|+\delta(t+2)}+\alpha\frac{d_t(v)}{|E_t|}.
\end{equation}
Therefore, by (\ref{1})
\begin{equation*}
    d_{t+1}^2(v_i)=d_t^2(v_i)+2d_t(v_i)\Delta d_{t+1}(v_i) +(\Delta d_{t+1}(v_i))^2,
\end{equation*}
and by (\ref{2})
\begin{multline}\label{EdsquareGt}
\E[d_{t+1}^2(v_i)\mid \mathcal{G}_t]=d_t^2(v_i)+2 d_t(v_i)\E[\Delta d_{t+1}(v_i)\mid\mathcal{G}_t] +\E[(\Delta d_{t+1}(v_i))^2\mid \mathcal{G}_t]=\\
d_t^2(v_i)+2\Big(\frac{(1-\alpha)d_t^2(v_i)}{2|E_t|+\delta(t+2)}+\frac{(1-\alpha)\delta d_t(v_i)}{2|E_t|+\delta(t+2)}+\frac{\alpha d_t^2(v_i)}{|E_t|}\Big) +\E[(\Delta d_{t+1}(v_i))^2\mid \mathcal{G}_t]
\end{multline}

Our next step is to show that replacing $|E_t|$ by its expectation in the RHS of \eqref{EdsquareGt} yields good approximations in expectation, in the following way:
\begin{equation}\label{59}
\E\Big[\frac{d_t^2(v_i)}{2|E_t|+\delta(t+2)}\Big]= \E\Big[\frac{d_t^2(v_i)}{2t(\alpha+1)+\delta t}\Big]\Big(1+O\Big(\sqrt{\frac{\ln t}{t}}\Big)\Big),
\end{equation}

\begin{equation}\label{60}
\E\Big[\frac{d_t(v_i)}{2|E_t|+\delta(t+2)}\Big]= \E\Big[\frac{d_t(v_i)}{2t(\alpha+1)+\delta t}\Big]\Big(1+O\Big(\sqrt{\frac{\ln t}{t}}\Big)\Big),
\end{equation}
and
\begin{equation}\label{60plus}
\E\Big[\frac{d_t^2(v_i)}{|E_t|}\Big]= \E\Big[\frac{d_t^2(v_i)}{t(\alpha+1)}\Big]\Big(1+O\Big(\sqrt{\frac{\ln t}{t}}\Big)\Big).
\end{equation}
Since these are so similar, we focus on \eqref{59}.

Consider the event $B_t$ defined by
\[B_t:=\{||E_t|-t(\alpha+1)-1|\leq\sqrt{t\ln t}\}.\]
Since $(|E_t|,t\geq0)$ is essentially a Bernoulli random walk, we have by Hoeffding's inequality, for any $x>0$,
\[\Prob(||E_t|-t(\alpha+1)-1|\geq x)\leq 2\,\mathrm{exp}\left( -\frac{2x^2}{t}\right).\]
Taking $x=\sqrt{t\ln t}$ yields $\Prob(B_t^c)\leq 2t^{-2}$.

Observe that since $t+1\leq|E_t|\leq 2t+1$ and $1\leq d_t(v_i)\leq 2|E_t|$, then $ct^{-1}\leq \frac{d_t^2(v_i)}{2|E_t|+\delta(t+2)}\leq Ct$, for some $c,C>0$ and $\tilde{c}t^{-1}\frac{d_t(v_i)}{2|E_t|+\delta(t+2)}\leq \tilde{C}$, for some $\tilde{c},\tilde{C}>0$.

Using this we have
\begin{equation}\label{53}
0\leq \E\Big[\frac{d_t^2(v_i)}{2|E_t|+\delta(t+2)} \mathbbm{1}_{\{|E_t|\notin B_t\}}\Big]\leq \frac{2Ct}{t^2}=\frac{2C}{t},
\end{equation} and
\begin{align}\label{54}
\E&\Big[\frac{d_t^2(v_i)\mathbbm{1}_{\{|E_t|\in B_t\}}}{2|E_t|+\delta(t+2)} \Big]\geq \E\Big[\frac{d_t^2(v_i)\mathbbm{1}_{\{|E_t|\in B_t\}}}{2(t(\alpha+1)+1+\sqrt{t\ln t})+\delta(t+2)}\Big]\nonumber\\
&= \E\Big[\frac{d_t^2(v_i)}{2(t(\alpha+1)+1+\sqrt{t\ln t})+\delta(t+2)}\Big] -\E\Big[\frac{d_t^2(v_i)\mathbbm{1}_{\{|E_t|\notin B_t\}}}{2(t(\alpha+1)+1+\sqrt{t\ln t})+\delta(t+2)}\Big]
\nonumber\\ &\geq  \E\Big[\frac{d_t^2(v_i)}{2(t(\alpha+1)+1+\sqrt{t\ln t})+\delta(t+2)}\Big] -\frac{2C}{t}.
\end{align}
Note that for $t$ large enough,
\begin{multline}\label{55}
\frac{\mathbbm{1}_{\{|E_t|\in B_t\}}}{2|E_t|+\delta(t+2)}\leq \frac{1}{2(t(\alpha+1)+1-\sqrt{t\ln t})+\delta(t+2)}\leq \frac{1}{2t(\alpha+1)+\delta t}\Big(1+C_2\sqrt{\frac{t}{\ln t}}\Big),
\end{multline} 
and
\begin{align}\label{56}
   \frac{1}{2(t(\alpha+1)+1+\sqrt{t\ln t})+\delta(t+2)}\geq \frac{1}{2t(\alpha+1)+\delta t}\Big(1-C_1\sqrt{\frac{t}{\ln t}}\Big)
\end{align}
where $C_2=4/[2(\alpha+1)+\delta-1]>0$, and $C_1=5/[2(\alpha+1)+\delta]>0$.

Then for $t$ large enough, by (\ref{53}) and (\ref{55})
\begin{align}\label{57}
\E\Big[\frac{d_t^2(v_i)}{2|E_t|+\delta(t+2)} \Big]&\leq
\E\Big[\frac{d_t^2(v_i)}{2t(\alpha+1)+\delta t}\Big]\Big(1+C_2\sqrt{\frac{\ln t}{t}}\Big)+ \frac{2C}{t}\nonumber\\
&\leq \E\Big[\frac{d_t^2(v_i)}{2t(\alpha+1)+\delta t}\Big]\Big(1+2C_2\sqrt{\frac{\ln t}{t}}\Big),
\end{align}
and by (\ref{53}), (\ref{54}) and (\ref{56})
\begin{align}\label{58}
\E\Big[\frac{d_t^2(v_i)}{2|E_t|+\delta(t+2)}\Big]&\geq
\E\Big[\frac{d_t^2(v_i)}{2t(\alpha+1)+\delta t}\Big]\Big(1-C_1\sqrt{\frac{\ln t}{t}}\Big) - \frac{C}{t}\nonumber\\
&\geq \E\Big[\frac{d_t^2(v_i)}{2t(\alpha+1)+\delta t}\Big]\Big(1-2C_1\sqrt{\frac{\ln t}{t}}\Big).
\end{align}
Hence we have \eqref{59}. \eqref{60} and \eqref{60plus} are obtained in analogous way.

Taking expectations in both sides of the equality in \eqref{EdsquareGt} and by (\ref{59}) and (\ref{60}), we obtain that $\E[d_{t+1}^2(v_i)]$ is equal to
\begin{equation}\label{Edsquare2}
\begin{aligned}
\E[d_t^2(v_i)]+
2\Bigg(\E\Big[\frac{(1-\alpha) d_t^2(v_i)}{2t(\alpha+1)+\delta t}\Big]+\E\Big[\frac{\delta(1-\alpha) d_t(v_i)}{2t(\alpha+1)+\delta t}\Big]+\E\Big[\frac{\alpha d_t^2(v_i)}{t(\alpha+1)}\Big]\Bigg)\Big(1+O\Big(\sqrt{\frac{\ln t}{t}}\Big)\Big) \\ +\E[(\Delta d_{t+1}(v_i))^2]\\
=\E[d_t^2(v_i)]\Bigg(1+\frac{2[(\alpha+1)^2+\alpha\delta]}{(\alpha+1)[2(\alpha+1)+\delta]t}+O\Big(\frac{\sqrt{\ln t}}{t\sqrt{t}}\Big)\Bigg)+
\frac{2\delta(1-\alpha)\E[d_t(v_i)]}{[2(\alpha+1)+\delta]t}\Big(1+O\Big(\sqrt{\frac{\ln t}{t}}\Big)\Big)\\+\E[(\Delta d_{t+1}(v_i))^2].
\end{aligned}
\end{equation}

Since $\Delta d_{t+1}(v_i)$ takes values in $\{0,1\}$, then $\E[\Delta d_{t+1}^2(v_i)\mid \mathcal{G}_t]$ is given by the RHS of \eqref{2}, and since $\sum_{i=1}^{|V_t|}d_t(v_i)=2|E_t|$, then $\sum_{i=1}^{|V_t|}\E[\Delta d_{t+1}^2(v_i)\mid \mathcal{G}_t]<2\alpha+1$, so $\sum_{i=1}^{|V_t|}\E[\Delta d_{t+1}^2(v_i)]<2\alpha+1$. Moreover, using that $\sum_{i=1}^{|V_t|}\E[d_t(v_i)]=2\E[|E_t|]=2[t(\alpha+1)+1]$ and by (\ref{Edsquare2}) we achieve
\begin{align*}
    \E[\tilde{C}_{t+1}]= \E[\tilde{C}_{t}]\Bigg(1+\frac{2[(\alpha+1)^2+\alpha\delta]}{(\alpha+1)[2(\alpha+1)+\delta]t}+O\Big(\frac{\sqrt{\ln t}}{t\sqrt{t}}\Big)\Bigg)+\Theta(1).
\end{align*}
Expanding this we obtain $\E[\tilde{C}_{t}]$ equals to
\begin{align}\label{CtildeExpan}
\prod_{s=1}^{t-1}\Bigg(1+\frac{2A}{s}+O\Big(\frac{\sqrt{\ln s}}{s\sqrt{s}}\Big)\Bigg)\E[\tilde{C}_1]+ 
\Theta(1)\sum_{s=1}^{t-1}\prod_{r=s}^{t-1}\Bigg(1+\frac{2A}{r}+O\Big(\frac{\sqrt{\ln r}}{r\sqrt{r}}\Big)\Bigg)
\end{align}
where $A=\frac{[(\alpha+1)^2+\alpha\delta]}{(\alpha+1)[2(\alpha+1)+\delta]}$.


Note that $\E[(\tilde{C}_1]=6(1+\alpha).$ Moreover, since for $s\geq 1$
$$
\prod_{r=s}^{t-1}\Bigg(1+\frac{2A}{r}+O\Big(\frac{\sqrt{\ln r}}{r\sqrt{r}}\Big)\Bigg)=\Theta\Big(\exp\Big(\sum_{r=s}^{t-1}\frac{2A}{r}\Big)\Big)=\Theta\Big(\frac{t}{s}\Big)^{2A},
$$
we have
\begin{align*}
\sum_{s=1}^{t-1}\prod_{r=s}^{t-1}\Bigg(1+\frac{2A}{r}+O\Big(\frac{\sqrt{\ln r}}{r\sqrt{r}}\Big)\Bigg)=
  \begin{cases}
   \Theta(t\ln t),& \text{ if } 2A=1\\
   \Theta(t),& \text{ if } 2A<1\\
   \Theta(t^{2A}),& \text{ if } 2A>1,
  \end{cases}
\end{align*}
and the same asymptotic holds for $\tilde{C}_t$:
\begin{align*}
\E[\tilde{C}_t]=
  \begin{cases}
   \Theta(t\ln t),& \text{ if } 2A=1\\
   \Theta(t),& \text{ if } 2A<1\\
   \Theta(t^{2A}),& \text{ if } 2A>1.
  \end{cases}
\end{align*}
Finally, since $C_t=\frac{1}{2}\tilde{C}_t-|E_t|$ and $\E[|E_t|]=t(\alpha+1)+1$, we achieve the desired result for $C_t$. 
\end{proof}

\noindent\textit{Proof of Theorem \ref{ThmExpC2}.}
We will start by obtaining an upper and a lower bound for $\E[\mathcal{C}_2(t)]$. To get the upper bound we will use 
that $T_t\leq t$, and $\E[T_t]=\alpha t$, so
\begin{align*}
  \E[\mathcal{C}_2(t)]=\E\Big[\frac{3T_t}{C_t}\Big]&\leq 3\E\Big[\frac{T_t}{\E[C_t]}\Big]+3\Big|\E\Big[\frac{T_t}{C_t}-\frac{T_t}{\E[C_t]}\Big] \Big| \\
  &=3\frac{\E[T_t]}{\E[C_t]}+3 \frac{t}{\E[C_t]}\Big|\E\Big[\frac{\E[C_t]-C_t}{C_t}\Big]\Big|\\
  &\leq 3\frac{\alpha t+t}{\E[C_t]}. 
\end{align*}
    For the lower bound of $\E[\mathcal{C}_2(t)]$, define the event $D_t=\{T_t:|T_t-\alpha t|\leq \sqrt{t\ln t}\}$. Since in $D_t$, $T_t\geq \alpha t-\sqrt{t\ln t}$, then
    \begin{align*}
        \E[\mathcal{C}_2(t)]&=\E\Big[3\frac{T_t}{C_t}\mathbbm{1}_{\{T_t\in D_t\}}\Big]+ \E\Big[3\frac{T_t}{C_t}\mathbbm{1}_{\{T_t\notin D_t\}}\Big]\\
        &\geq 
        \E\Big[3\frac{T_t}{C_t}\mathbbm{1}_{\{T_t\in D_t\}}\Big]\\
        &\geq 3\E\Big[\frac{\alpha t-\sqrt{t\ln t}}{C_t}\Big]\\
        &= \frac{3(\alpha t-\sqrt{t\ln t})}{\E[C_t]}.
    \end{align*}
The result of this theorem follows by Lemma \ref{lemEPt}.    
\qed


\vskip0.5cm

\section{Conclusions}
It remains a challenge to find a model which features a power-law degree distribution with $2<\gamma<3$ and still features positive global clustering, which is what is observed in several real-world networks. In \cite{Ostroumova} for example, the authors say that the reason for this
is that the out-degree (i.e. number of parents) is a fixed constant value, and consequently the number of triangles grows too slowly. This applies "morally" to our model, as the out-degree here is random, but has fixed expectation $1+\alpha$. For this reason, we believe that introducing more complex structures at each instant of time, for example a 3 or 4 dimensional simplex would only moderately improve the clustering. It is also mentioned in~\cite{Ostroumova} that a way to overcome this obstacle is to increase the out-degree. In \cite{CooperPralat11}, the authors studied a preferential attachment model with out-degree given by a function of $t$, $f(t)=t^c$. They proved that  for $0<c < 1$, the degree sequence of the process exhibits a power-law with parameter $\gamma = (3 -c)/(1 -c)>3$, though the clustering was not studied. Obtaining models featuring positive global clustering with $\gamma<3$ would require moving far away enough from the framework of \cite{Ostroumova}.

\vskip0.5cm
\section*{Acknowledgements} The authors would like to thank the Heilbronn Institute for Mathematical Research (HIMR) and the UKRI/EPSRC 
Additional Funding Programme for Mathematical Sciences, for their support under the Small Grant Call.


\end{document}